\documentclass[onecolumn,draftclsnofoot]{IEEEtran}
\usepackage{cite}
\usepackage{amsmath,amssymb,amsfonts}
\usepackage{algorithmic}
\usepackage{graphicx}
\usepackage{textcomp}
\usepackage[caption=false]{subfig}
\usepackage[flushleft]{threeparttable}
\def\BibTeX{{\rm B\kern-.05em{\sc i\kern-.025em b}\kern-.08em
    T\kern-.1667em\lower.7ex\hbox{E}\kern-.125emX}}    
\usepackage{xcolor}

\makeatletter
\usepackage{algorithm}
\newcommand\fs@norules{\def\@fs@cfont{\mathbfseries}\let\@fs@capt\floatc@ruled
	\def\@fs@pre{}%
	\def\@fs@post{}%
	\def\@fs@mid{\kern3pt}%
	\let\@fs@iftopcapt\iftrue}
\makeatother
\floatstyle{norules}
\restylefloat{algorithm}
\allowdisplaybreaks

\begin{document}
\title{\fontsize{14}{16}\selectfont
A Memory-efficient Implementation of Perfectly Matched Layer with Smoothly-varying Coefficients in Discontinuous Galerkin Time-Domain Method}
\author{Liang Chen, Mehmet Burak Ozakin, Shehab Ahmed, and Hakan Bagci
\thanks{
	The authors are with the Division of Computer, Electrical, and Mathematical Science and Engineering, King Abdullah University of Science and Technology (KAUST), Thuwal 23955-6900, Saudi Arabia (e-mails:\{liang.chen, mehmet.ozakin, shehab.ahmed, hakan.bagci\}@kaust.edu.sa).
}
\thanks{This work is supported by the King Abdullah University of Science and Technology (KAUST) Office of Sponsored Research (OSR) under Award No 2019-CRG8-4056. The authors would like to thank the KAUST Supercomputing Laboratory (KSL) for providing the required computational resources.}}
\maketitle

\begin{abstract}
Wrapping a computation domain with a perfectly matched layer (PML) is one of the most effective methods of imitating/approximating the radiation boundary condition in Maxwell and wave equation solvers. Many PML implementations often use a smoothly-increasing attenuation coefficient to increase the absorption for a given layer thickness, and, at the same time, to reduce the numerical reflection from the interface between the computation domain and the PML. In discontinuous Galerkin time-domain (DGTD) methods, using a PML coefficient that varies within a mesh element requires a different mass matrix to be stored for every element and therefore significantly increases the memory footprint. In this work, this bottleneck is addressed by applying a weight-adjusted approximation to these mass matrices. The resulting DGTD scheme has the same advantages as the scheme that stores individual mass matrices, namely higher accuracy (due to reduced numerical reflection) and increased meshing flexibility (since the PML does not have to be defined layer by layer) but it requires significantly less memory.
\end{abstract}



\section{Introduction}
\IEEEPARstart{P}{erfectly} matched layer (PML)~\cite{Berenger1994, Chew1994} is often used in finite difference~\cite{Taflove2005}, finite element~\cite{Jin2015FEM}, and discontinuous Galerkin time-domain (DGTD) methods~\cite{Hesthaven2002, Cockburn2004, Lu2004, Hesthaven2008, Gedney2009, Cohen2017} to imitate/approximate the radiation boundary condition (i.e., truncate a unbounded physical domain to a finite computation domain) while solving Maxwell equations or the wave equation. 


The performance of the PML depends on the attenuation coefficient (which is implemented as conductivity in Maxwell equations) and the thickness of the layer. Increasing either one or both of them increases the absorption inside the PML. However, in practice, one cannot use a constant and high conductivity as it would increase the numerical reflection at the interface between the PML and computation domain, or use a very thick layer since it would increase the computational cost. Therefore, a smoothly increasing conductivity profile is often used to achieve both high absorption and small numerical reflection~\cite{Berenger1994, Chew1994, Chew1996, Taflove2005, Jin2015FEM, Berenger2007book, Gedney2011}.

\begin{figure}[!b]
	\centerline{\includegraphics[width=0.79\columnwidth]{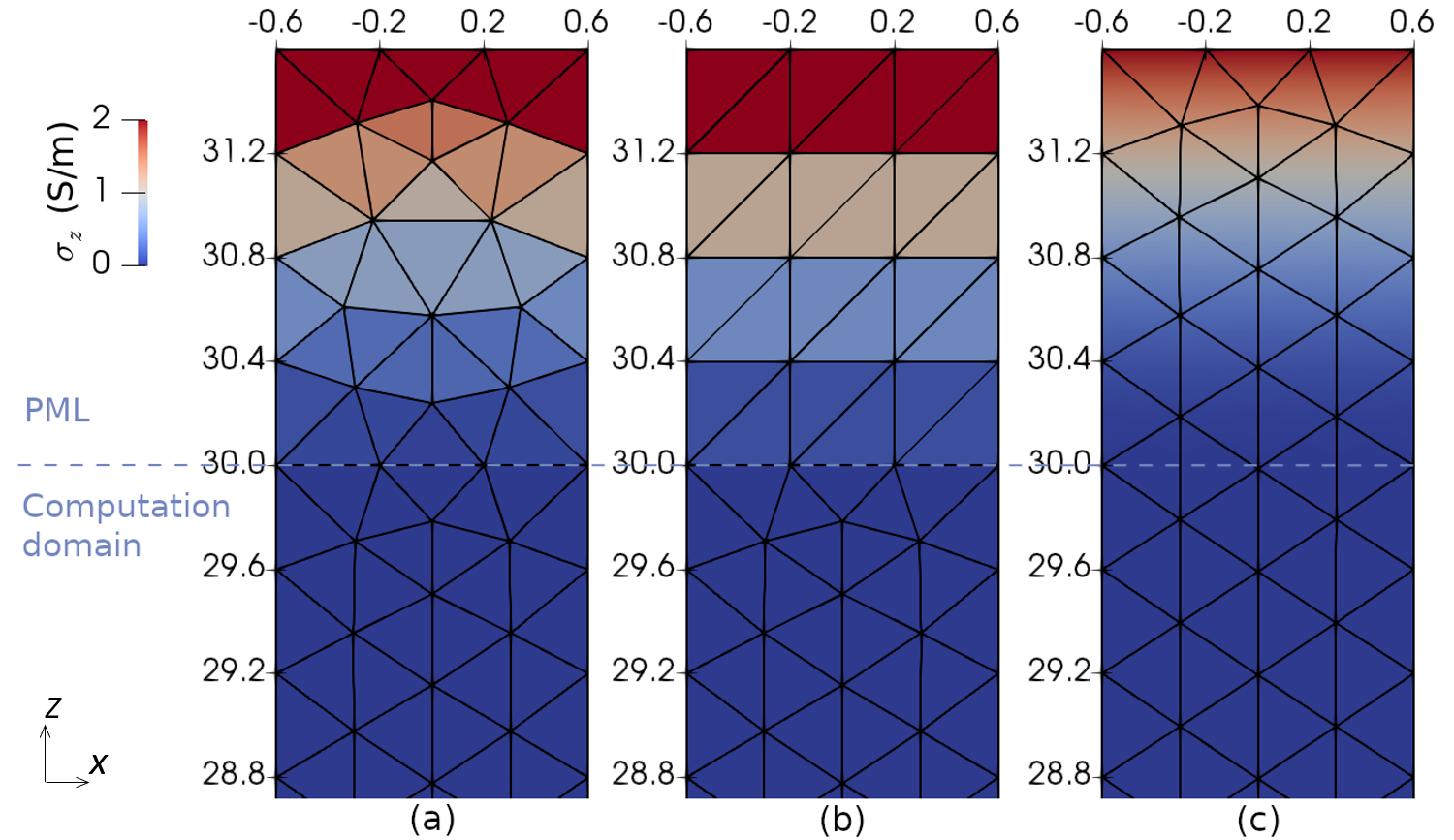}}
	\caption{Implementation of the PML with smoothly-increasing conductivity. (a) Conductivity is constant in elements (paved mesh). (b) Conductivity is constant in elements (layered mesh). (c) Conductivity is allowed to vary in elements (paved mesh). Mesh does not have to conform to the interface between the PML and the computation domain.}
	\label{Profile}
\end{figure}

In DGTD, the PML conductivity profile can be implemented in two different ways. The first method assumes that the conductivity is constant in a given element. This constant might for example be set to the value of the conductivity profile at the center of the element. Implementation of this method is rather straightforward since the mass matrices of different elements only differ by a constant (for linear elements)~\cite{Hesthaven2008}. However, the conductivity profile becomes discontinuous between neighboring elements [Fig.~\ref{Profile} (a)] and the element surfaces that are not parallel to the PML interface lead to large reflections and destroy the high-order accuracy of the solution. One workaround is to build layered (tetrahedral) meshes [Fig.~\ref{Profile} (b)] or use orthogonal (hybrid) meshes inside the PML~\cite{Lu2004, Chen2018}, and accordingly set a layered conductivity profile. But this makes the setups of the computation domain and the PML rather tedious since one needs to control the mesh and the conductivity on all face, edge, and corner regions of the PML. Moreover, to reduce the numerical reflection by decreasing the conductivity discontinuity between neighboring mesh (or conductivity) layers, their number has to be increased.

The second method allows the conductivity to vary inside a given element (following the increasing conductivity profile inside the PML). It should be noted here that (higher-order) DGTD allows for sampling of the material properties at the sub-elemental level. In this case, the behavior of the PML is determined only by the conductivity samples within the elements and therefore the mesh interfaces can be aligned arbitrarily [Fig.~\ref{Profile} (c)] without adversely affecting its performance. {\color{black}Indeed, this second approach can improve the PML performance~\cite{Angulo2015, Angulo2014, Sankaran2007thesis}.} However, the conductivity varying within the elements results in an element-dependent mass matrix. Therefore, a direct implementation requires a different mass matrix (or its inverse) to be stored for every element~\cite{Gedney2009, Lu2004, Niegemann2009, Angulo2015, Angulo2014, Chen2018}. This significantly increases DGTD’s memory footprint. For example, for the stretched-coordinate (SC)-PML~\cite{Gedney2009}, the memory required to store the mass matrices scales with $15K_{\mathrm{PML}} \times N_p^2$, where ${N_p}$ is the number of interpolating nodes in each element, $K_{\mathrm{PML}}$ is the number of elements in the PML, and $15$ comes from the five material-dependent coefficients in the update equations and three Cartesian components of the vector field. In contrast, for the first method, where the conductivity is assumed constant in a given element, this memory requirement scales with $K_{\mathrm{PML}}$ since only the constant conductivity of each element and a single reference mass matrix are stored. 

In this work, a memory-efficient method to implement the SC-PML with smoothly-varying attenuation coefficients in DGTD is developed. The proposed method allows the conductivity to vary inside the elements and constructs the resulting local mass matrices using a weight-adjusted approximation (WAA)~\cite{Chan2017}. Compared to the direct implementation that is briefly described above, the proposed method reduces the memory requirement to $15K_{\mathrm{PML}} \times N_q$, where $N_q\sim N_p$, while maintaining the PML performance. It should be noted here that the WAA has been proved to be energy-stable and preserve the high-order convergence of DGTD~\cite{Chan2017, Guo2020, Shukla2020}. Indeed, numerical examples presented here also show the proposed method maintains the higher-order accuracy of the solution. Additionally, it is demonstrated that the PML with smoothly-increasing conductivity profile as implemented with the proposed method performs better than the PML implemented using element-wise constant conductivity profile.

\section{Formulation}
\subsection{WAA-DGTD for SC-PML}
Maxwell equations in stretched-coordinates for a source-free and lossless medium can be expressed as~\cite{Chew1994}
\begin{align}
\label{MaxwellSCE}
- j\omega \mu {\mathbf{H}} = {\nabla _e} \times {\mathbf{E}}\\
\label{MaxwellSCH}
j\omega \epsilon {\mathbf{E}} = {\nabla _h} \times {\mathbf{H}}
\end{align}
where $\mathbf{E}$ and $\mathbf{H}$ are the electric and magnetic fields, $\epsilon$ and $\mu$ are the permittivity and the permeability, $\omega$ is the frequency, and
\begin{equation}
{\nabla _e} = {\nabla _h} = \hat x\frac{1}{{{s_x}}}\frac{\partial }{{\partial x}} + \hat y\frac{1}{{{s_y}}}\frac{\partial }{{\partial y}} + \hat z\frac{1}{{{s_z}}}\frac{\partial }{{\partial z}}.
\end{equation}
The coordinate-stretching variables ${s_{u}}$, $u \in \{x,y,z\}$, are defined as~\cite{Chew1994, Gedney2009, Gedney2011}
\begin{equation}
\label{su}
{s_u (u)} = {\kappa _u (u)} + \frac{{{\sigma _u (u)}}}{{j\omega {\varepsilon _0}}}
\end{equation}
where $\kappa _u$ and $\sigma _u$ are one-dimensional positive real functions along direction $u$. Here, $\sigma _u$ is the attenuation coefficient that ensures absorption inside the SC-PML and $\kappa _u$ changes the propagation speed {\color{black} (and also increases the absorption for evanescent waves)}. It is well known that using smoothly-increasing $\sigma _u$ and $\kappa _u$ reduces the numerical reflection from the interface between the PML and the computation domain while maintaining a high absorption rate, and therefore improves the PML performance~\cite{Chew1996, Berenger2007book, Gedney2011}.

The time-domain update equations in SC-CPML can be expressed as~\cite{Gedney2009}
\begin{align}
\label{upH}
& {\partial _t}\ddot a \cdot \mu {\mathbf{H}} = - \nabla \times {\mathbf{E}} - \ddot b \cdot \mu {\mathbf{H}} - \ddot c \cdot \mu {{\mathbf{P}}^H} \\
\label{upE}
& {\partial _t}\ddot a \cdot \varepsilon {\mathbf{E}} = \nabla \times {\mathbf{H}} - \ddot b \cdot \varepsilon {\mathbf{E}} - \ddot c \cdot \mu {{\mathbf{P}}^E} \\
\label{upPH}
& {\partial _t}{{\mathbf{P}}^H} = {\ddot \kappa ^{ - 1}}{\mathbf{H}} - \ddot d{{\mathbf{P}}^H} \\
\label{upPE}
& {\partial _t}{{\mathbf{P}}^E} = {\ddot \kappa ^{ - 1}}{\mathbf{E}} - \ddot d{{\mathbf{P}}^E}
\end{align}
where ${{\mathbf{P}}^E}$ and ${{\mathbf{P}}^H}$ are auxiliary variables introduced to avoid computationally costly temporal convolutions while converting~\eqref{MaxwellSCE} and~\eqref{MaxwellSCH} into time domain~\cite{Gedney2009} and $\ddot a$, $\ddot b$, $\ddot c$, $\ddot d$, and $\ddot \kappa$ are diagonal tensors with entries defined as
\begin{align}
\nonumber & {a_{uu}} = \frac{{{\kappa _v}{\kappa _w}}}{{{\kappa _u}}}, \quad
 {b_{uu}} = \frac{1}{{{\kappa _u}{\varepsilon _0}}}({\sigma _v}{\kappa _w} + {\sigma _w}{\kappa _v} - {a_{uu}}{\sigma _u})\\
\label{tensor} & {c_{uu}} = \frac{{{\sigma _v}{\sigma _w}}}{{\varepsilon _0^2}} - {b_{uu}}\frac{{{\sigma _u}}}{{{\varepsilon _0}}}, \quad
 {d_{uu}} = \frac{{{\sigma _u}}}{{{\kappa _u}{\varepsilon _0}}}, \quad
 {\kappa _{uu}} = {\kappa _u}.
\end{align}
In~\eqref{tensor} and the rest of the text $(u,v,w)$ follows the permutation $(x,y,z)$ $\rightarrow$ $(y,z,x)$ $\rightarrow(z,x,y)$.

Following the standard nodal discontinuous Galerkin method~\cite{Hesthaven2008}, the computation domain and the PML are discretized into $K$ elements with volumetric support $\Omega_k$ and boundary surface $\partial \Omega_k$, and in each element ${\mathbf{E}}$, ${\mathbf{H}}$, ${{\mathbf{P}}^E}$ and ${{\mathbf{P}}^H}$ are expanded using the Lagrange polynomials $\ell _i({\mathbf{r}})$~\cite{Hesthaven2008}, $i=1,\cdots,N_p$, where ${N_p} = (p + 1)(p + 2)(p + 3)/6$ is the number of interpolating nodes and $p$ is the order of the Lagrange polynomials. Finally, Galerkin testing yields the semi-discrete system of equations as
\begin{align}
\label{disH} 
& {\partial _t}{{\bar{H}}_k} = - ({\bar{M}}_k^a)^{-1} [{\bar{M}}_k^b{\bar{H}}_k + {\bar{M}}_k^c{{\bar{P}}_k^H}
 + {\mu_k^{-1}} \bar{\mathbb{C}}_k
({{\bar{E}}_k},{{\bar{E}}_{k'}},{{\bar{H}}_k},{{\bar{H}}_{k'}})] \\
\label{disE} 
& {\partial _t}{{\bar{E}}_k} = - ({\bar{M}}_k^a)^{-1} [{\bar{M}}_k^b{\bar{E}}_k + {\bar{M}}_k^c{{\bar{P}}_k^E}
 - {\varepsilon_k^{-1}} \bar{\mathbb{C}}_k({{\bar{H}}_k},{{\bar{H}}_{k'}}, {{\bar{E}}_k},{{\bar{E}}_{k'}})] \\
\label{disPH}
& {\partial _t}{\bar{P}}_k^H = {{\bar{M}}_k^{-1}} ({\bar{M}}_k^{{1/\kappa}}{\bar{H}}_k - {\bar{M}}_k^d{{\bar{P}}_k^H}) \\
\label{disPE}
& {\partial _t}{\bar{P}}_k^E = {{\bar{M}}_k^{-1}} ({\bar{M}}_k^{{1/\kappa}}{\bar{E}}_k - {\bar{M}}_k^d{{\bar{P}}_k^E}).
\end{align}
Here, $\bar{H}_k$, $\bar{E}_k$, $\bar{P}_k^H$, and $\bar{P}_k^E$ are vectors storing the unknown time-dependent coefficients of the relevant basis functions, $\bar{M}_k$ and $\bar{M}_k^{\alpha}$, $\alpha \in \{a,b,c,d,{1/\kappa}\}$, are the mass matrices with entries
\begin{align}
\label{mass0}
& {\bar{M}_k}(i,j) = \int_{{\Omega _k}} {{\ell _i}({\mathbf{r}}){\ell _j}({\mathbf{r}})} d{\mathbf{r}} \\
\label{mass}
& {\bar{M}_k^{\alpha,u}}(i,j) = \int_{{\Omega _k}} {\alpha_{uu}(\mathbf{r}){\ell _i}({\mathbf{r}}){\ell _j}({\mathbf{r}})} d{\mathbf{r}}
\end{align}
$\bar{\mathbb{C}}_k({f_k},{f_{k'}},{g_k},{g_{k'}})$ denotes the curl operator with its component along direction $u$ is given by
\begin{equation*}
\bar{\mathbb{C}}_k^u({f_k},{f_{k'}},{g_k},{g_{k'}}) = {{\bar{S}}_k^v}{f_k^w} - {{\bar{S}}_k^w}{f_k^v} + {{\bar{F}}_{k}}{\mathbb{F}^u}({f_k},{f_{k'}},{g_k},{g_{k'}})
\end{equation*}
where $u \in \{ x,y,z\}$, $(f,g)\in\{(\bar{E}, \bar{H}), (\bar{H}, \bar{E})\}$, $\mathbb{F}^u$ is the component of the numerical flux along direction $u$, which in general involves unknowns from the current element $k$ and its neighboring element $k'$~\cite{Hesthaven2008, Chen2020steadystate, Chen2019multiphysics}, $\bar{S}_k$ and $\bar{F}_k$ are the stiffness and the face mass matrices with entries
\begin{align}
\label{stiff}
& \bar{S}_k^u (i,j) = \int_{{\Omega _k}} {{\ell _i}({\mathbf{r}})\frac{{d{\ell _j}({\mathbf{r}})}}{{du }}} d{\mathbf{r}} \\
\label{lift}
& {\bar{F}_k}(i,j) = \oint_{\partial {\Omega _k}} {{\ell _i}({\mathbf{r}}){\ell _j}({\mathbf{r}})d{\mathbf{r}}}
\end{align}
and $\epsilon _k$ and $\mu _k$ are the permittivity and the permeability (assumed constant) in each element. Note that, the nodal DG framework~\cite{Hesthaven2008} is used here, but the proposed method can be easily extended to vector DG methods~\cite{Cockburn2004, Lu2004, Gedney2009, Li2015IBC, Li2017dispersive, Li2018graphene}.

For linear elements, the mass matrices ${\bar{M}_k}$ in~\eqref{mass0} are simply scaled versions of the mass matrix ${\bar{M}}$ defined on the reference element, ${\bar{M}_k = J_k {\bar{M}}}$, where $J_k$ is the Jacobian of the coordinate transformation between element $k$ and the reference element. Hence, only ${\bar{M}}$ and (scalar constant) $J_k$ are stored. Similarly, in~\eqref{mass}, if $\alpha_{uu}(\mathbf{r})$ is assumed constant inside the elements, i.e., $\alpha_{uu}(\mathbf{r}) = \alpha_{uu}^k$, then ${\bar{M}_k^{\alpha,u}}=\alpha_{uu}^k{\bar{M}_k}=\alpha_{uu}^k J_k {\bar{M}}$ and one only needs to store $\alpha_{uu}^k$ in addition to ${\bar{M}}$ and $J_k$. In this case, \eqref{disH}-\eqref{disPE} can be implemented as efficiently as the case without the PML~\cite{Hesthaven2008, Liu2012, Chen2019discontinuous}.

However, if $\alpha_{uu}(\mathbf{r})$ is allowed to vary inside the elements, ${\bar{M}_k^{\alpha,u}}$ are different in different elements and in general there is no simple relationship between these different mass matrices. One has to store every one of these mass matrices (or their inverse). As an alternative, the mass matrix can be recomputed at each time step, but this would significantly increase the cost of time marching~\cite{Hesthaven2008}.
The memory required to store the mass matrices ${\bar{M}_k^{\alpha,u}}$ in~\eqref{disH}-\eqref{disPE} scales with $3 \times 5\times N_p^2$ per element, where $3$ is the number of the $(x,y,z)$ components of the vector field, $5$ is the number of the coefficients $a(\mathbf{r})$, $b(\mathbf{r})$, $c(\mathbf{r})$, $d(\mathbf{r})$, and $\kappa(\mathbf{r})$. Note that this memory requirement is significantly higher than that of storing the unknowns coefficients of the basis functions, which scales with $12 \times N_p$ in the PML.

To reduce the memory requirement of implementing~\eqref{mass} with $\alpha_{uu}(\mathbf{r})$ allowed to vary inside the elements, WAA~\cite{Chan2017} is used. It has been shown that with this approximation DGTD retains provable energy-stability and high-order accuracy~\cite{Chan2017, Guo2020, Shukla2020}. Note that in the above SC-PML formulation, directly multiplying~\eqref{upH} and \eqref{upE} with $\ddot{a}^{-1}$ on both sides reduces the number of element-dependent mass matrices to $4$. But this would result in a non-conservative form, whose solution is neither provably energy-stable nor provably  high-order accurate~\cite{Hesthaven2008, Chan2017}.

First, a weight-adjusted inner product is introduced to approximate the parameter-weighted inner product in the expression of the mass matrix~\cite{Chan2017}. The mass matrix, which is associated with the element $k$ and a locally varying coefficient $\alpha_k(\mathbf{r})$, is approximated as
\begin{align}
\label{WAmass}
{\bar{M}_k^{\alpha}} \approx \bar{M}_k (\bar{M}_k^{1/\alpha})^{-1} \bar{M}_k.
\end{align}
Since $({\bar{M}_k^{\alpha}})^{-1}$ is used in~\eqref{disH}-\eqref{disPE} (for $\alpha=a$), one needs to  calculate $\bar{M}_k^{1/\alpha}$. Under the nodal DG framework~\cite{Hesthaven2008},
\begin{align}
\nonumber
\bar{M}_k^{1/\alpha}(i,j) & = J_k \int_{{\Omega _k}} {\alpha_k^{-1}(\mathbf{r}){\ell _i}({\mathbf{r}}){\ell _j}({\mathbf{r}})} d{\mathbf{r}} \\
& \approx J_k \sum_q \ell_i(\mathbf{r}_q) {w_q}{\alpha_k^{-1}(\mathbf{r}_q)}  \ell_j(\mathbf{r}_q)
\end{align}
where $\mathbf{r}_q$, $q=1, ..., N_q$, are the Gaussian quadrature nodes corresponding to the quadrature rules of degree $2p+1$ and $w_q$ are the corresponding weights. Hence,
\begin{align}
\label{WAmass1}
\bar{M}_k^{1/\alpha} = J_k \bar{V}_q^T \bar{w}_q \bar{\alpha}_k^{-1} \bar{V}_q
\end{align}
where $\bar{V}_q$ is an interpolation matrix defined on the reference element, $\bar{V}_q = \bar{V}_I \bar{V}^{-1}$, $\bar{V}_I$ and $\bar{V}$ are generalized Vandermonde matrices with entries $\bar{V}_I(q,i)=\phi_i(\mathbf{r}_q)$ and $\bar{V}(j,i)=\phi_i(\mathbf{r}_j)$, respectively, $\phi_i(\mathbf{r})$ is the $i$-th orthonormal polynomial basis~\cite{Hesthaven2008}, $\bar{w}_q=\mathrm{diag} \{w_1, ..., w_{N_q}\}$ is also element-independent, and $\bar{\alpha}_k = \mathrm{diag} \{ \alpha_k(\mathbf{r}_1), ... , \alpha_k(\mathbf{r}_{N_q}) \}$ is a diagonal matrix containing the coefficients evaluated at the quadrature nodes.

Inverting~\eqref{WAmass} and substituting~\eqref{WAmass1} yield
\begin{align}
\nonumber
({\bar{M}_k^{\alpha}})^{-1} & \approx \bar{M}_k^{-1} \bar{M}_k^{1/\alpha} \bar{M}_k^{-1}\\
\nonumber
& = {\bar{M}^{-1}} \bar{V}_q^T \bar{w}_q \bar{\alpha}_k^{-1} \bar{V}_q \bar{M}_k^{-1} \\
\label{WAmassInv}
& = \bar{P}_q \bar{\alpha}_k^{-1} \bar{V}_q \bar{M}_k^{-1}.
\end{align}
Here, $\bar{P}_q = \bar{M}^{-1}\bar{V}_q^T \bar{w}_q$ is introduced to simplify the implementation. In~\eqref{WAmassInv}, $\bar{P}_q$ and $\bar{V}_q$ are defined on the reference element, and ${\bar{M}_k^{-1}}$ is a scaled version of the reference matrix, ${\bar{M}_k^{-1}}= J_k^{-1} {\bar{M}^{-1}}$.

The update equations~\eqref{disH}-\eqref{disPE} contain multiplications between element-dependent mass matrices. To reduce the number of arithmetic operations, the following operators are defined
\begin{align}
\label{Mb}
& \tilde{M}_k^{b} = {({\bar{M}}_k^{a})^{-1} {\bar{M}}_k^{b}} = \bar{P}_q \bar{a}_k^{-1} \bar{V}_q \bar{P}_q \bar{b}_k \bar{V}_q \\
\label{Mc}
& \tilde{M}_k^{c} = {({\bar{M}}_k^{a})^{-1} {\bar{M}}_k^{c}} = \bar{P}_q \bar{a}_k^{-1} \bar{V}_q \bar{P}_q \bar{c}_k \bar{V}_q \\
\label{Md}
& \tilde{M}_k^{d} = {{\bar{M}}_k^{-1} {\bar{M}}_k^{d}} = \bar{P}_q \bar{d}_k \bar{V}_q\\
\label{Mkappa}
& \tilde{M}_k^{1/\kappa} = {{\bar{M}}_k^{-1} {\bar{M}}_k^{1/\kappa}} = \bar{P}_q \bar{\kappa}_k^{-1} \bar{V}_q
\end{align}
where~\eqref{WAmassInv} is used for $\alpha=a$ and~\eqref{WAmass1} is used for $\alpha \in \{b, c, d, {1/\kappa}\}$. These operators can be directly used on the right hand sides of~\eqref{disH}-\eqref{disPE}. Substituting~\eqref{Mb}-\eqref{Mkappa} into~\eqref{disH}-\eqref{disPE} yields
\begin{align}
\label{disH3}
& {\partial _t}{{\bar{H}}_k} =  - \bar{P}_q \bar{a}_k^{-1} \bar{V}_q [\bar{P}_q (\bar{b}_k \bar{V}_q \bar{H}_k + \bar{c}_k \bar{V}_q \bar{P}_k^H )
+ {\mu_k^{-1}} \bar{\mathbb{C}}_k ({{\bar{E}}_k},{{\bar{E}}_{k'}},{{\bar{H}}_k},{{\bar{H}}_{k'}}) ] \\
\label{disE3}
& {\partial _t}{{\bar{E}}_k} =  - \bar{P}_q \bar{a}_k^{-1} \bar{V}_q [\bar{P}_q (\bar{b}_k \bar{V}_q \bar{E}_k + \bar{c}_k \bar{V}_q \bar{P}_k^E )
- {\epsilon_k^{-1}} \bar{\mathbb{C}}_k ({{\bar{H}}_k},{{\bar{H}}_{k'}},{{\bar{E}}_k},{{\bar{E}}_{k'}}) ] \\
\label{disPH3}
& {\partial _t}{\bar{P}}_k^H = \bar{P}_q (\bar{\kappa}_k^{-1} \bar{V}_q {\bar{H}}_k - \bar{d}_k \bar{V}_q {\bar{P}}_k^H ) \\
\label{disPE3}
& {\partial _t}{\bar{P}}_k^E = \bar{P}_q (\bar{\kappa}_k^{-1} \bar{V}_q {\bar{E}}_k - \bar{d}_k \bar{V}_q {\bar{P}}_k^E ).
\end{align}
Equations~\eqref{disH3}-\eqref{disPE3} can be implemented in a matrix-free manner just like it is done in classical DG implementations~\cite{Hesthaven2008, Gedney2009, Cohen2017, Liu2012, Sirenko2012, Chen2019discontinuous, Chen2019unitcell}. 


\subsection{Computational complexity}
\label{complexity}
In DGTD with explicit time marching, all operations are localized within the elements. The memory required to store the mass matrices in the direct implementation of~\eqref{disH}-\eqref{disPE} scales with $K_{\mathrm{PML}} \times 15 N_p^2$, where $15$ comes from the number of unknown components times the number of different mass matrices associated with different coefficients and $K_{\mathrm{PML}}$ is the number of elements in the PML. In the WAA formulation~\eqref{disH3}-\eqref{disPE3}, the memory requirement reduces to $(K_{\mathrm{PML}} \times 15 N_q) + 2 N_p N_q$, where $15 N_q$ comes from the number of unknown components times the number of coefficient samples at the quadrature points and $2 N_p N_q$ comes from $\bar{V}_q$ and $\bar{P}_q$ defined on the reference element. For simplicial quadrature rules that are exact for up to polynomials of degree $2p+1$, $N_q\! \sim \! N_p$~\cite{Cools1999, Xiao2010}.

To compare the number of arithmetic operations required by the two implementations, one should first note that the curl operator $\bar{\mathbb{C}}$ is the same in both formulations. Computation of $\bar{\mathbb{C}}$ requires those of the spatial derivatives and the numerical flux~\cite{Chen2020float, Sirenko2018, Chen2020hybridizable}. Here, the memory access time is much more significant than the time required to carry out these computations because data from neighboring elements, which are discontinuous in memory, is required. Therefore, only the times required to complete the arithmetic operations of the remaining terms are compared. For the same reason, in practice, the time required to compute $\bar{\mathbb{C}}$ dominates the overall time required by the time marching, and the difference in the numbers of arithmetic operations as estimated below for the remaining terms is less significant (see the example in Section~\ref{Examples}).

In~\eqref{disH}, the three matrix-vector multiplications and two vector-vector additions require $3 N_p^2$ multiplication operations and $2 N_p$ addition operations, respectively. In~\eqref{disH3}, the multiplication of $\bar{V}_q$ with a vector of length $N_p$, and the multiplication of $\bar{P}_q$ with a vector of length $N_q$ require $N_q N_p$ multiplication operations. The multiplication of a diagonal matrix with a vector (such as $\bar{b}_k \bar{v}$) requires $N_q$ multiplication operations. As a result, excluding the computation of $\bar{\mathbb{C}}$,~\eqref{disH3} requires $5 N_q N_p + 3 N_q$ multiplications and $N_q + N_p$ additions. For the auxiliary variable, the cost of~\eqref{disPH} is $3 N_p^2$ multiplications and $N_p$ subtractions, while~\eqref{disPH3} requires $3 N_q N_p + 2 N_q$ multiplications and $N_q$ subtractions. One can see the number of operations in the WAA implementation is slightly higher than that in the direct implementation. But as mentioned above the time required by these operations is smaller than the time required to compute $\bar{\mathbb{C}}$, and therefore overall times required by the two implementations are not that different.

\section{Numerical Examples}\label{Examples}
In this section, the accuracy and the efficiency of the proposed WAA formulation are compared to those of the traditional PML implementations using numerical examples. To this end, four PML configurations/implementations are considered in these examples: (i) $\sigma_u$ and/or $\kappa_u$, $u \in \{x,y,z\}$, are assumed constant inside the elements on a paved mesh (EC-paved) [Fig.~\ref{Profile}(a)], (ii) $\sigma_u$ and/or $\kappa_u$, $u \in \{x,y,z\}$, are assumed constant inside the elements on a layered mesh (EC-layered) [Fig.~\ref{Profile}(b)], (iii) $\sigma_u$ and/or $\kappa_u$, $u \in \{x,y,z\}$, are allowed to vary inside the elements on a paved mesh (SV-paved) [Fig.~\ref{Profile}(c)], and (iv) same configuration in (iii) but implemented using the proposed method with the WAA (SV-WAA-paved) [Fig.~\ref{Profile}(c)]. In all implementations, the order of the Lagrange polynomials $p \in \{1, 2, 3, 4, 5,\}$, which results in $N_p \in \{4, 10, 20, 35, 56\}$. 

For configuration (i), the constant values in a given element are obtained by sampling $\sigma_u$ and $\kappa_u$, $u \in \{x,y,z\}$, at that element's node that is farthest away from the PML interface (along the $\pm u$-direction). For configuration (ii), to ensure that the element surfaces are strictly parallel to the axes, the PML mesh is built layer by layer and constant values in a given layer are obtained by sampling $\sigma_u$ and $\kappa_u$, $u \in \{x,y,z\}$, at the outermost surface of that layer (along the $\pm u$-direction). For the WAA in implementation (iv), the order of the Gaussian quadrature rule is $2p$, resulting in $N_q \in \{4, 11, 23, 44, 74\}$~\cite{Xiao2010}.

In all examples, the background medium is free space and the excitation is a plane wave with electric field $\mathbf{E}(z,t)=E_0\mathbf{\hat{x}}G(t-z/c_0)$, where $E_0=1 \mathrm{V/m}$, $c_0$ is the speed of light in free space, and $G(t) = e^{(t-t_0)^2/4\tau^2}$ is a base-band Gaussian pulse with $\tau=66.67\,\mathrm{ps}$ and $t_0=15\tau$. The average edge lengths of all meshes used under this excitation are $0.4\,\mathrm{cm}$. 

First, the reflection of a plane wave normally incident on the PML is computed. The computation domain is a rectangular box with dimensions $1.2\,\mathrm{cm} \times 1.2\,\mathrm{cm} \times 60\, \mathrm{cm}$. Perfect electric conductor (PEC) and periodic boundary conditions are used on the outer boundary of the PML that is located perpendicular to the $z$ direction and on the computation domain boundaries perpendicular to the $x$ and $y$ directions, respectively. The plane wave excitation is introduced on surface $z=0$ and propagates in the $+z$-direction.
The domain is long enough to ensure that the reflected field is well-separated from the incident one, and therefore the reflection from the PML is simply measured by the peak value of the reflected field’s amplitude. The conductivity profile is described by $\sigma_z(z)=\sigma_{\mathrm{max}}[(z-z_0)/L_{z}]^{p_{\sigma}}$, where $z_0$ is the $z$-coordinate on the interface between PML and the computation domain, $L_z$ is the thickness of the PML and $ p_{\sigma}$ is the order of the profile. Note that $\sigma_z(z)$ is nonzero only when $|z|>|z_0|$. The values of these parameters are $z_0=\pm 30$ $\mathrm{cm}$, $L_{z}=1.6$ $\mathrm{cm}$, and $p_{\sigma}=1$, and also $\kappa_z(z)=1$ both inside the PML and the computation domain.

\begin{figure}[!t]
	\centerline{\includegraphics[width=0.79\columnwidth]{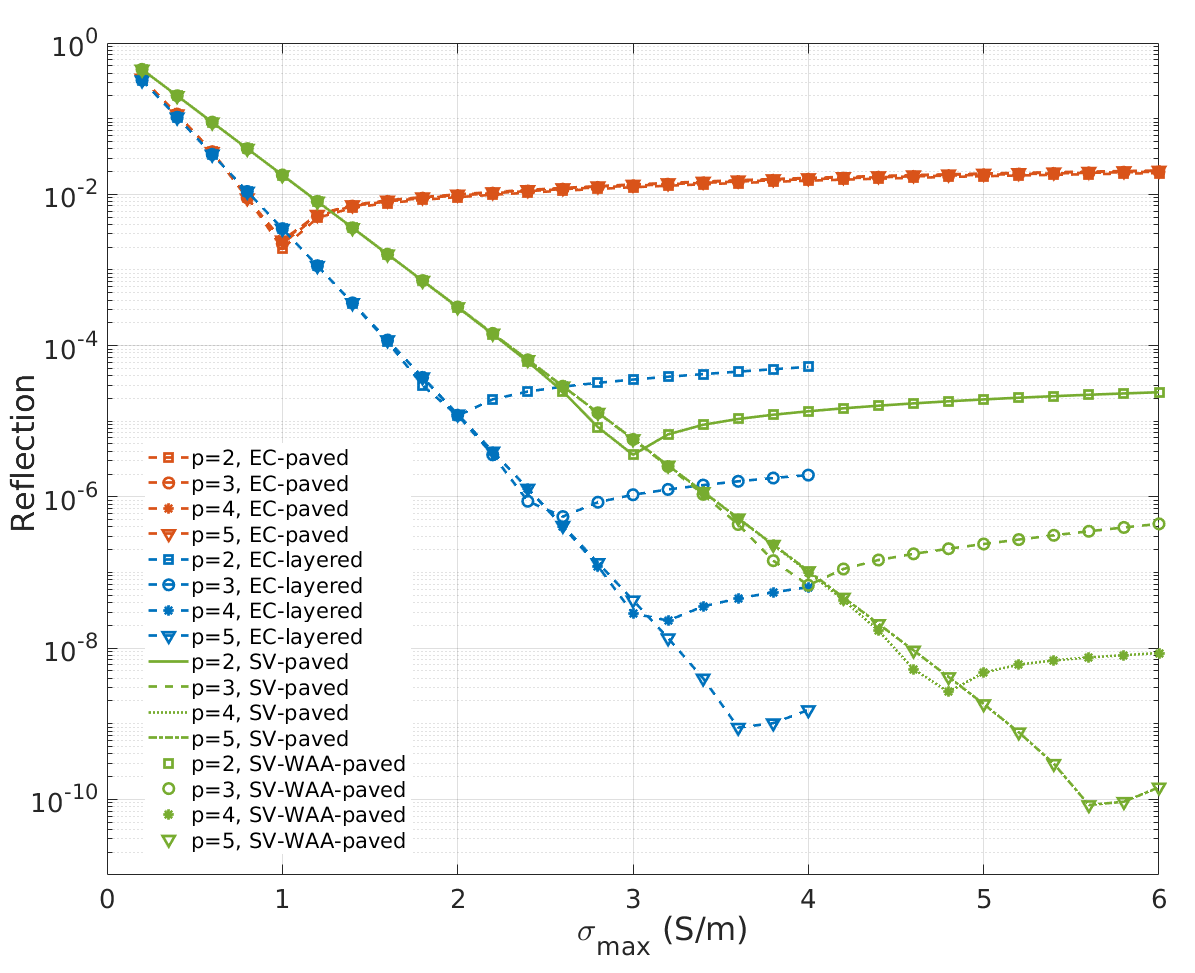}}
	\caption{Peak value of the reflected field’s amplitude versus $\sigma_{\mathrm{max}}$ for different PML configurations/implementations and orders of Lagrange polynomials ($p$). The plane wave  is normally incident on the PML.}
	\label{PW}
\end{figure}

In this example, four configurations/implementations are considered: EC-paved, EC-layered, SV-paved, and SV-WAA-paved. Their performances are compared for $p \in \{2, 3, 4, 5\}$. For all four groups of simulations, $\sigma_{\mathrm{max}}$ is scanned to find the minimum reflection that can be obtained for each case. Fig.~\ref{PW} shows that with increasing $\sigma_{\mathrm{max}}$, the reflection first decreases exponentially and then increases gradually. This is observed for all configurations/implementations and all values of $p$. When $\sigma_{\mathrm{max}}$ is small, the overall reflection is dominated by the reflection from the PEC boundary simply because the absorption inside the PML is not high enough. Therefore, in this regime, increasing $\sigma_{\mathrm{max}}$ elevates the absorption and reduces the amplitude of the wave reflected back into the computation domain exponentially. The numerical reflection (which is smaller than the reflection from the PEC boundary for small $\sigma_{\mathrm{max}}$) increases with increasing $\sigma_{\mathrm{max}}$~\cite{Chew1996}, and starts dominating the overall reflection as demonstrated in the figure by the gradual increase after the minimum point.

For the EC-paved configuration, the reflection stays at a high level and does not decrease with increasing $p$. This is because of the large reflections from unoriented internal element surfaces. For the EC-layered, SV-layered, and SV-WAA-paved configurations, high-order convergence is observed, i.e., the reflection keeps on decreasing exponentially with increasing $p$. Still, the reflection for the SV-paved and SV-WAA-paved configurations is about $15~\mathrm{dB}$ smaller than that for the EC-layered configuration. Note that this higher accuracy comes with the ease of meshing since a layered mesh (and conductivity profile) is not needed. Finally, Fig.~\ref{PW} also shows that the SV-WAA-paved implementation performs exactly the same as the SV-paved direct implementation, which verifies the accuracy of the proposed method. 

\begin{figure}[!t]
	\centering
	\subfloat[\label{Article_SphereL}]{\includegraphics[width=0.499\columnwidth]{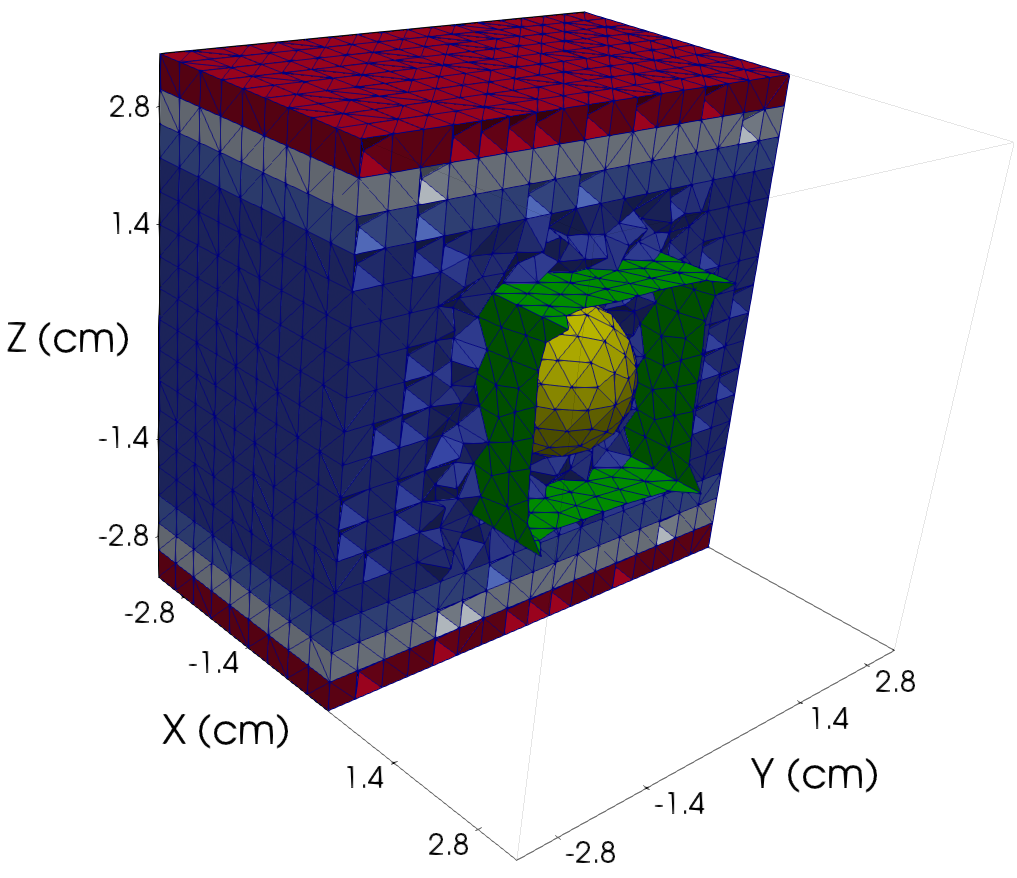}}
	\subfloat[\label{Article_SphereR}]{\includegraphics[width=0.499\columnwidth]{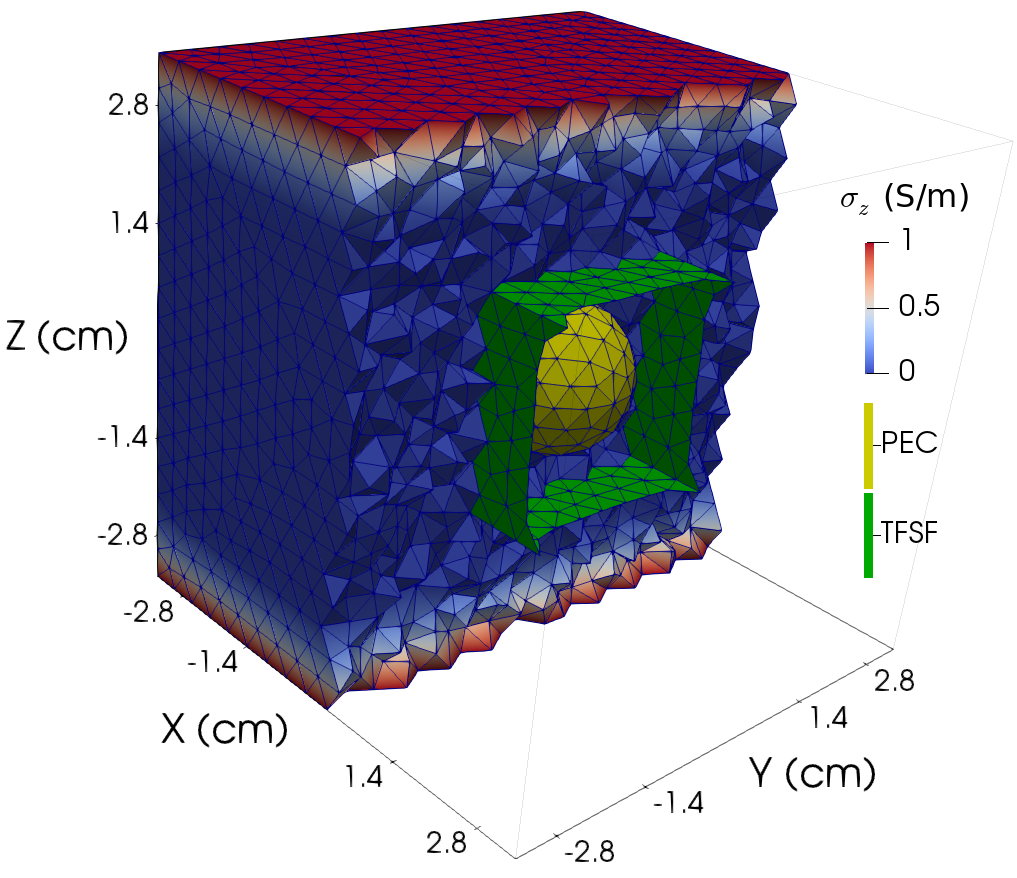}}
	\caption{ Computation domains, meshes, and PML conductivity profiles (represented by color) used for the (a) EC-layered and (b) SV-paved and SV-WAA-paved configurations.}
	\label{Article_Sphere}
\end{figure}

Next, scattering from a PEC sphere of radius $1$ $\mathrm{cm}$ is considered. The computation domains and the PMLs for the EC-layered and SV-paved and SV-WAA-paved configurations are shown in Figs.~\ref{Article_Sphere} (a) and (b), respectively. The plane wave excitation is introduced on the total-field scattered-field (TFSF) surface [shown in green in Figs.~\ref{Article_Sphere} (a) and (b)]. The conductivity function is $\sigma_u(u)=\sigma_{\mathrm{max}}[(u-u_0)/L_{u}]^{p_{\sigma}}$, $u \in \{x,y,z\}$, $u_0 \in \{x_0, y_0, z_0\}$, where $u_0$ is the $u$-coordinate on the interface between PML and the computation domain, $L_u$ is the thickness of the PML along the $\pm u$ direction, and $ p_{\sigma}$ is the order of the profile. Note that $\sigma_u(u)$ is nonzero only when $|u|>|u_0|$. The values of these parameters are $x_0=y_0=z_0=\pm 2.2$ $\mathrm{cm}$, and $L_x=L_y=L_z=1.2$ $\mathrm{cm}$.

Because the distance between the sphere surface and the PML is short, possibly-evanescent scattered waves enter the PML with high grazing angles. A varying $\kappa_u$, $u \in \{x,y,z\}$, profile is employed to help with the absorption of these evanescent waves~\cite{Berenger2007book, Gedney2011}: $\kappa_u(u)=1+(\kappa_{\mathrm{max}}-1)[(u-u_0)/L_{u}]^{p_{\sigma}}$ with $\kappa_{\mathrm{max}}=2$. Note that inside the computation domain, $\kappa_u(u)=1$.  In this example, using the PEC or the first-order absorbing boundary condition~\cite{Angulo2015} on the outer boundary of the PML gives similar results. The results presented here are obtained with the PEC boundary condition. 

Three configurations/implementations are considered here: EC-layered, SV-paved, and SV-WAA-paved. For the EC-layered configuration, the thickness of each PML mesh layer is $0.4$ $\mathrm{cm}$ [Fig.~\ref{Article_Sphere} (a)]. Note that, for this example, generation of these layers is rather tedious since in the corner region one has to align all layer/element surfaces in all three directions. In contrast, for the SV-paved and SV-WAA-paved configurations [same mesh is used -- Fig.~\ref{Article_Sphere} (b)], $\sigma_u$ and $\kappa_u$ values are simply obtained by sampling the corresponding profile functions at the nodes of the elements. This significantly simplifies the setups of the computation domain and the PML since even an explicit interface between the computation domain and the PML is not required [see Fig.~\ref{Article_Sphere} (b)]. The performances of the three configurations are compared for $p \in \{1, 3, 4\}$ and $ p_{\sigma} \in \{1, 2\}$. 

\begin{figure}[!t]
	\centerline{\includegraphics[width=0.79\columnwidth]{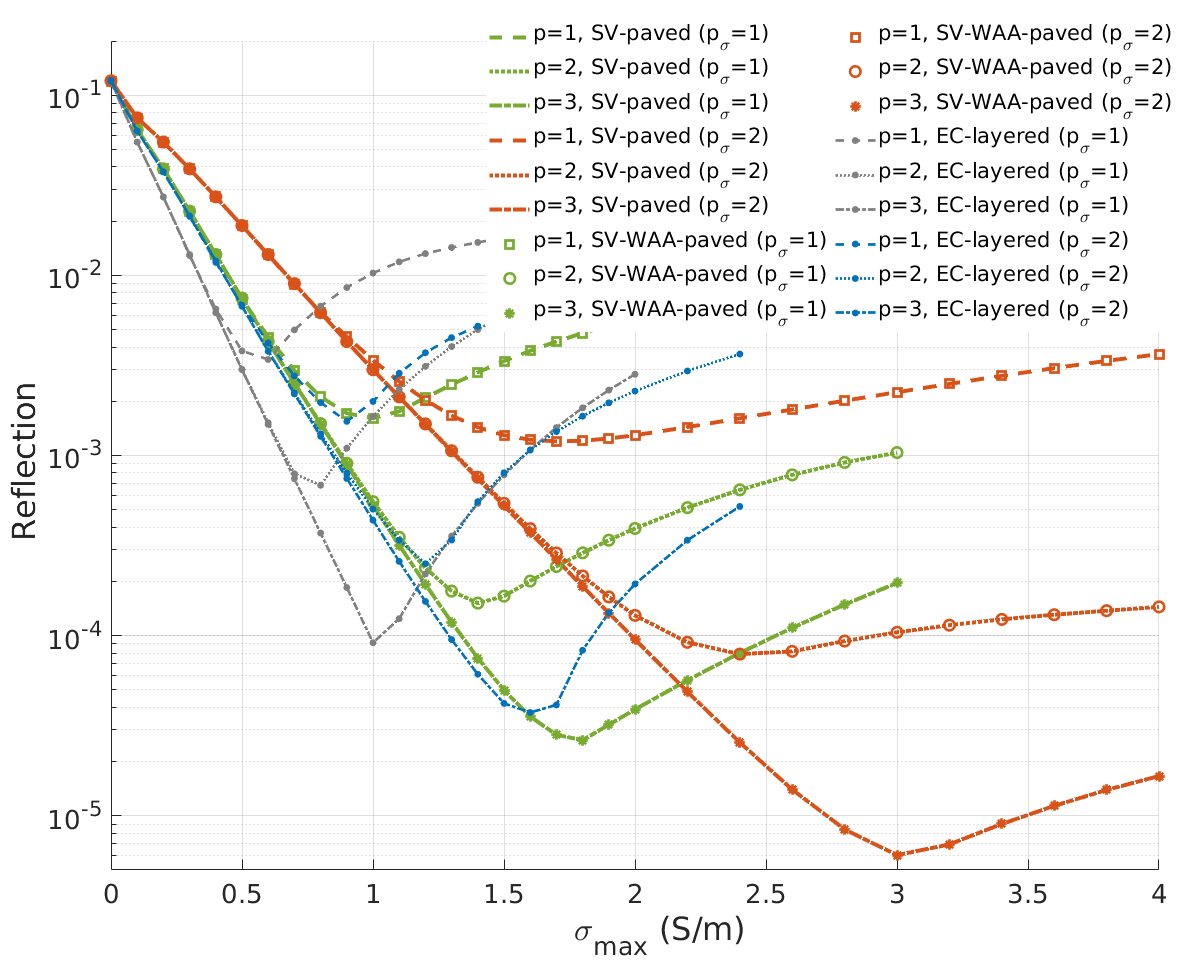}}
	\caption{Reflection of the scattered field versus $\sigma_{\mathrm{max}}$ for different PML configurations/implementations, orders of the Lagrange polynomials ($p$), and orders of the PML profile ($p_{\sigma}$). The scatterer is a PEC sphere of radius $1\, \mathrm{cm}$. }
	\label{SphereReflectionSigMax}
\end{figure}

$\sigma_{\mathrm{max}}$ is scanned to find the minimum reflection that could be reached for each case. Note that in this example ``reflection'' is defined as the peak value of the absolute difference between the fields computed at a probe point for the above cases and corresponding reference fields computed at the same point. $10$ different probe points (placed either in the TF or in the SF region) have been tested and the results are consistent for all of them. The results below correspond to the probe point at $(1.0 \,\mathrm{cm}, 1.0\, \mathrm{cm}, 1.0\, \mathrm{cm})$. These reference fields are computed under the same excitation but the distance between the sphere surface and the PML is extended to $12$ $\mathrm{cm}$. To ensure that the discretization errors are at the same level, the meshes in the overlapped regions between the actual computation domains and the extended ones are kept exactly the same, the average edge lengths of the meshes in the extended region are kept the same as those in the actual computation domains, and the solutions are obtained using the same value of $p$.

Fig.~\ref{SphereReflectionSigMax} plots the reflection for the three cases with different values of $p$ and $p_{\sigma}$ versus $\sigma_{\mathrm{max}}$. Clearly, the SV-paved and SV-WAA-paved configurations perform better than the EC-layered configuration for every value of $p_{\sigma}$. The best performance is obtained with $p_{\sigma}=2$. Note that further increasing $p_{\sigma}$ degrades the PML performance for all configurations/implementations and all values of $p$ since high conductivity values only appear at the very end of the PML when $p_{\sigma}$ is high. 
Fig.~\ref{SphereReflectionSigMax} also shows that the SV-WAA-paved implementation performs exactly the same as the SV-paved direct implementation, which means the error caused by the WAA of the mass matrices is below the level of the discretization error.

\begin{table}[!t]
	\centering
	\begin{threeparttable}
		\renewcommand{\arraystretch}{1.1}
		\centering
		\caption{Computational costs of the SV-paved and SV-WAA-paved implementations for  different orders of the Lagrange polynomials ($p$) \tnote{*}.}
		\label{cost}
		\setlength{\tabcolsep}{3pt}
		\begin{tabular}{ c | c | c | c | c | c | c }
			\hline 
			$p$ & $N_p$ & $N_q$ & \multicolumn{2}{|c|}{memory (KB)} & \multicolumn{2}{|c}{CPU time per step (s)} \\ 
			\hline
			&    & & SV-paved     & SV-WAA-paved & SV-paved   & SV-WAA-paved\\ \hline
			1 & 4  & 4  & 378,660    & 267,928   & 1.652716 & 2.341525 \\ \hline
			2 & 10 & 11 & 1,274,424  & 498,508   & 4.083981 & 5.960960 \\ \hline
			3 & 20 & 23 & 4,126,640  & 894,936   & 9.606642 & 15.73330 \\ \hline
			4 & 35 & 44 & 11,583,440 & 1,513,140 & 19.64900 & 30.16986 \\ \hline
			5 & 56 & 74 & 28,410,000 & 2,291,608 & 78.56877 & 105.1317 \\ \hline
		\end{tabular}
		\smallskip
		\scriptsize
		\begin{tablenotes}
			\item[*] {Tested on a workstation with Intel Xeon(R) E5-2680 v4 CPU and 128GB memory. A single process is used. $K=72,762$ and $K_{\mathrm{PML}}=52,657$.}
		\end{tablenotes}
	\end{threeparttable}
\end{table}
Table.~\ref{cost} compares the computational cost of the SV-paved and SV-WAA-paved implementations. With increasing $p$, the memory requirement increases dramatically for the SV-paved direct implementation but only modestly for the SV-WAA implementation. For $p=5$, the memory requirement of the SV-paved implementation is $12.4$ times that of the SV-WAA-paved implementation. The computation time required by the SV-WAA-paved implementation per time step is slightly larger than that required by the SV-paved implementation due to the increased number of arithmetic operations (see Section.~\ref{complexity}). It should be also be noted here that, in practice, a DGTD algorithm is usually parallelized. The difference in times required for updating different elements is relatively small and can be easily compensated by allocating a smaller number of elements for those MPI processes containing PML elements. In the numerical results presented here, assigning a weight of $2$ for PML elements in ParMetis~\cite{metis1998, Chen2019parallel} yields a good load-balance.

\section{Conclusion}
A PML implementation that allows the attenuation coefficient to vary inside the discretization elements yields a smaller numerical reflection from the interface between the PML and the computation domain and significantly simplifies the meshing process. However, these advantages come at the cost of increased memory footprint since a different mass matrix has to be stored for every discretization element. In this work, this memory requirement is reduced by applying WAA to the mass matrices without abandoning the advantages listed above. Indeed, numerical results demonstrate that the PML with smoothly-increasing conductivity profile as implemented with the proposed method performs better than the PML implemented using element-wise constant conductivity profile and that the higher-order accuracy of the solution is maintained.

The proposed method is especially useful for simulations running on shared-memory systems where the high memory requirement of smoothly-varying PMLs could be a bottleneck. For simulations running on distributed-memory systems, the memory requirement of a single computing node is also reduced and a better load-balance could be reached with a slightly adjusted weight in the domain partition.





\end{document}